\documentclass{amsart}
\usepackage{amssymb}
\newtheorem{theorem}{Theorem}
\newtheorem{lemma}[theorem]{Lemma}

\newtheorem{corollary}[theorem]{Corollary}
\newtheorem{definition}[theorem]{Definition}
\newcommand{\Map}{\mathrm{Map}}
\newcommand{\CP}{\mathbb{CP}}
\newcommand{\F}{\mathbb{F}}
\title{Fiber sums of genus 2 Lefschetz fibrations}
\author{D. Auroux}
\address{Department of Mathematics, M.I.T., Cambridge MA 02139, USA}
\email{auroux@math.mit.edu}

\begin{document}
\begin{abstract}
Using the recent results of Siebert and Tian about the holomorphicity
of genus 2 Lefschetz fibrations with irreducible singular fibers, we show
that any genus 2 Lefschetz fibration becomes holomorphic after fiber
sum with a holomorphic fibration.
\end{abstract}

\maketitle

\section{Introduction}
Symplectic Lefschetz fibrations have been the focus of a lot of attention
since the proof by Donaldson that, after blow-ups, every compact symplectic
manifold admits such structures \cite{Do}. Genus 2 Lefschetz fibrations,
where the first non-trivial topological phenomena arise, have been
particularly studied. Most importantly, it has recently been shown by Siebert
and Tian that every genus 2 Lefschetz fibration without reducible fibers and
with ``transitive monodromy'' is holomorphic \cite{ST2}. The statement
becomes false if reducible singular fibers are allowed, as evidenced by
the construction by Ozbagci and Stipsicz \cite{OS} of genus 2 Lefschetz
fibrations with non-complex total space (similar examples have
also been constructed by Ivan Smith).

It has been conjectured by Siebert and Tian that any genus 2 Lefschetz
fibration should become holomorphic after fiber sum with sufficiently
many copies of the rational genus 2 Lefschetz fibration with 20 irreducible
singular fibers. The purpose of this paper is to
prove this conjecture by providing a classification of genus 2 Lefschetz
fibrations up to stabilization by such fiber sums. The result is the
following (see \S 2 and Definition 4 for notations):

\begin{theorem}
Let $F$ be any factorization of the identity element as a product of
positive Dehn twists in the mapping class group $\Map_2$. Then there
exist integers $\epsilon\in\{0,1\}$, $k\ge 0$ and $m\ge 0$ such that, 
for any large enough integer $n$, the factorization $F\cdot (W_0)^n$
is Hurwitz equivalent to $(W_0)^{n+k} \cdot (W_1)^\epsilon \cdot (W_2)^m$.
\end{theorem}

\begin{corollary}
Let $f:X\to S^2$ be a genus $2$ Lefschetz fibration. Then the fiber sum of
$f$ with sufficiently many copies of the rational genus $2$ Lefschetz
fibration with $20$ irreducible singular fibers is isomorphic to a 
holomorphic fibration.
\end{corollary}

\section{Mapping class group factorizations}

Recall that a Lefschetz fibration $f:X\to S^2$ is a fibration admitting
only isolated singularities, all lying in distinct fibers of $f$, and
near which a local model for $f$ in orientation-preserving complex
coordinates is given by $(z_1,z_2)\mapsto z_1^2+z_2^2$. We will only
consider the case $\dim X=4$, where the smooth fibers are compact surfaces
(of genus $g=2$ in our case), and the singular fibers present nodal
singularities obtained by collapsing a simple closed loop (the {\it
vanishing cycle}) in the smooth fiber. The monodromy of the fibration
around a singular fiber is given by a positive Dehn twist along the
vanishing cycle.

Denoting by $q_1,\dots,q_r\in S^2$ the images of the singular fibers and
choosing a reference point in $S^2$, we can characterize the fibration $f$
by its {\it monodromy} $\psi:\pi_1(S^2-\{q_1,\dots,q_r\})\to \Map_g$, where
$\Map_g=\pi_0\mathrm{Diff}^+(\Sigma_g)$ is the mapping class group of a
genus $g$ surface. It is a classical result (cf.\ \cite{Kas}) that the
monodromy morphism $\psi$ is uniquely determined up to conjugation by an
element of $\Map_g$ and a braid acting on $\pi_1(S^2-\{q_i\})$, and 
that it determines the isomorphism class of the Lefschetz fibration $f$.

While all positive Dehn twists along non-separating curves are mutually
conjugate in $\Map_g$, there are different types of twists along separating
curves, according to the genus of each component delimited by the curve.
When $g=2$, only two cases can occur: either the curve splits the surface
into two genus
$1$ components, or it is homotopically trivial and the corresponding
singular fiber contains a sphere component of square $-1$. The
latter case can always be avoided by blowing down the total space of the
fibration; if the blown-down fibration can be shown to be
holomorphic, then by performing the converse blow-up procedure we conclude
that the original fibration was also holomorphic. Therefore, in all the
following we can assume that our Lefschetz fibrations are {\it relatively
minimal}, i.e.\ have no homotopically trivial vanishing cycles.

The monodromy of a Lefschetz fibration can be encoded in a {\it mapping
class group factorization} by choosing an ordered system of generating loops
$\gamma_1,\dots,\gamma_r$ for $\pi_1(S^2-\{q_i\})$, such that each loop 
$\gamma_i$ encircles only one of the points $q_i$ and $\prod \gamma_i$
is homotopically trivial. The monodromy of the fibration along each of
the loops $\gamma_i$ is a Dehn twist $\tau_i$; we can then describe the
fibration in terms of the relation $\tau_1\cdot \ldots\cdot \tau_r=1$ in
$\Map_2$. The choice of the loops $\gamma_i$ (and therefore of the twists
$\tau_i$) is of course not unique, but any two choices differ by a sequence
of {\it Hurwitz moves} exchanging consecutive factors:
$\tau_i\cdot \tau_{i+1} \to (\tau_{i+1})_{\tau_i^{-1}}\cdot
\tau_i$ or $\tau_i\cdot \tau_{i+1} \to \tau_{i+1}\cdot
(\tau_i)_{\tau_{i+1}}$, where we use the notation
$(\tau)_{\phi}=\phi^{-1}\tau \phi$, i.e.\ if $\tau$ is a Dehn twist
along a loop $\delta$ then $(\tau)_\phi$ is the Dehn twist along the
loop $\phi(\delta)$.

\begin{definition}
A {\em factorization} $F=\tau_1\cdot\ldots\cdot \tau_r$ in $\Map_g$ is an
ordered tuple of positive Dehn twists. We say that two factorizations are 
{\em Hurwitz equivalent} $(F\sim F')$ if they can be obtained from each 
other by a sequence of Hurwitz moves.
\end{definition}

It is well-known that a Lefschetz fibration is characterized by a
factorization of the identity element in $\Map_g$, uniquely determined up 
to Hurwitz equivalence and simultaneous
conjugation of all factors by a same element of $\Map_g$.

\begin{figure}[b]
\begin{center}
\setlength{\unitlength}{1cm}
\begin{picture}(5,2)(0,-1)
\qbezier[200](0,0)(0,1.2)(1.5,1)
\qbezier[200](1.5,1)(2.5,0.85)(3.5,1)
\qbezier[200](3.5,1)(5,1.2)(5,0)
\qbezier[200](0,0)(0,-1.2)(1.5,-1)
\qbezier[200](1.5,-1)(2.5,-0.85)(3.5,-1)
\qbezier[200](3.5,-1)(5,-1.2)(5,0)
\qbezier[60](1,0)(1,0.3)(1.5,0.3)
\qbezier[60](2,0)(2,0.3)(1.5,0.3)
\qbezier[60](1,0)(1,-0.3)(1.5,-0.3)
\qbezier[60](2,0)(2,-0.3)(1.5,-0.3)
\qbezier[60](3,0)(3,0.3)(3.5,0.3)
\qbezier[60](4,0)(4,0.3)(3.5,0.3)
\qbezier[60](3,0)(3,-0.3)(3.5,-0.3)
\qbezier[60](4,0)(4,-0.3)(3.5,-0.3)
\put(1.5,0){\circle{1.1}}
\put(3.5,0){\circle{1.1}}
\qbezier[60](0,0)(0,-0.1)(0.5,-0.1)
\qbezier[60](1,0)(1,-0.1)(0.5,-0.1)
\qbezier[20](0,0)(0,0.1)(0.5,0.1)
\qbezier[20](1,0)(1,0.1)(0.5,0.1)
\qbezier[60](2,0)(2,-0.1)(2.5,-0.1)
\qbezier[60](3,0)(3,-0.1)(2.5,-0.1)
\qbezier[20](2,0)(2,0.1)(2.5,0.1)
\qbezier[20](3,0)(3,0.1)(2.5,0.1)
\qbezier[60](4,0)(4,-0.1)(4.5,-0.1)
\qbezier[60](5,0)(5,-0.1)(4.5,-0.1)
\qbezier[20](4,0)(4,0.1)(4.5,0.1)
\qbezier[20](5,0)(5,0.1)(4.5,0.1)
\qbezier[100](2.5,0.91)(2.6,0.91)(2.6,0)
\qbezier[100](2.5,-0.91)(2.6,-0.91)(2.6,0)
\qbezier[30](2.5,0.91)(2.4,0.91)(2.4,0)
\qbezier[30](2.5,-0.91)(2.4,-0.91)(2.4,0)
\put(-0.35,0){$\zeta_1$}
\put(0.8,0.5){$\zeta_2$}
\put(2.1,-0.35){$\zeta_3$}
\put(3.9,0.5){$\zeta_4$}
\put(5.1,0){$\zeta_5$}
\put(2.65,0.7){$\sigma$}
\end{picture}

\caption{}
\end{center}
\end{figure}

Let $\zeta_i$ ($1\le i\le 5$) and $\sigma$ be the Dehn twists represented
in Figure 1. It is well-known (cf.\ e.g.\ \cite{Bi}, Theorem
4.8) that $\Map_2$ admits the following presentation:

\begin{tabular}{ll}
-- generators: & $\zeta_1,\dots,\zeta_5$.\\
-- relations: & $\zeta_i\zeta_j=\zeta_j\zeta_i$ if $|i-j|\ge 2$;
$\zeta_i\zeta_{i+1}\zeta_i=\zeta_{i+1}\zeta_i\zeta_{i+1}$; \\
& $(\zeta_1\zeta_2\zeta_3\zeta_4\zeta_5)^6=1$;
$I=\zeta_1\zeta_2\zeta_3\zeta_4\zeta_5^2\zeta_4\zeta_3\zeta_2\zeta_1$ is
central; $I^2=1$.
\end{tabular}

It is easy to check that $\sigma$ can be expressed in terms of
the generators $\zeta_1,\dots,\zeta_5$ as $\sigma=(\zeta_1\zeta_2)^6=
(\zeta_4\zeta_5)^6=(\zeta_1\zeta_2)^3(\zeta_4\zeta_5)^3 I$.

We can fix a hyperelliptic structure on the genus 2 surface $\Sigma$, i.e.\ a
double covering map $\Sigma\to S^2$ (with 6 branch points), in such a way
that $\zeta_1,\dots,\zeta_5$ become the lifts of standard half-twists
exchanging consecutive branch points in $S^2$. The element $I$ then
corresponds to the hyperelliptic involution (i.e.\ the non-trivial
automorphism of the double covering). The fact that $I$ is central means
that every diffeomorphism of $\Sigma$ is compatible with the hyperelliptic
structure, up to isotopy. In fact, $\Map_2$ is closely related to the braid
group $B_6(S^2)$ acting on the branch points of the double covering.
The group $B_6(S^2)$
admits the following presentation  (cf.\ \cite{Bi}, Theorem 1.1):

\begin{tabular}{ll}
-- generators: & $x_1,\dots,x_5$ (half-twists exchanging two consecutive 
points).\\

-- relations: & $x_ix_j=x_jx_i$ if $|i-j|\ge 2$; $x_ix_{i+1}x_i=x_{i+1}x_i
x_{i+1}$; \\
& $x_1x_2x_3x_4x_5^2x_4x_3x_2x_1=1$.
\end{tabular}
\medskip

Consider a $S^2$-bundle $\pi:P\to S^2$, and a smooth curve $B\subset P$
intersecting a generic fiber in $6$ points, everywhere transverse to
the fibers of $\pi$ except for isolated nondegenerate complex tangencies.
The curve $B$ can be characterized by its {\it braid monodromy}, or 
equivalently by a factorization in the braid group $B_6(S^2)$, with
each factor a positive half-twist, defined by considering the motion 
of the $6$ intersection points of $B$ with the fiber of $\pi$
upon moving around the image of a tangency point. As before, this
factorization is only defined up to Hurwitz equivalence and simultaneous
conjugation (see also \cite{Moi} for the case of plane curves).

There exists a lifting morphism from $B_6(S^2)$ to $\Map_2/\langle
I\rangle$, defined by $x_i\mapsto \zeta_i$. Given a half-twist in $B_6(S^2)$,
exactly one of its two possible lifts to $\Map_2$ is a Dehn twist about a
non-separating curve. This allows us to lift the braid factorization
associated to the curve $B\subset P$ to a mapping class group factorization;
the product of the resulting factors is equal to $1$ if the homology class
represented by $B$ is divisible by two, and to $I$ otherwise. In the first
case, we can construct a genus 2 Lefschetz fibration by considering
the double covering of $P$ branched along $B$, and its monodromy is exactly
the lift of the braid monodromy of the curve $B$. This construction always
yields Lefschetz fibrations without reducible singular fibers; however,
if we additionally allow some blow-up and blow-down operations (on $P$ and
its double covering respectively), then we can also handle the case of
reducible singular fibers (see \S 3 below and \cite{ST1}). It is worth
mentioning that Siebert and Tian have shown the converse result:
given any genus 2 Lefschetz fibration, it can be realized as a double
covering of a $S^2$-bundle over $S^2$ (with additional blow-up and blow-down
operations in the case of reducible singular fibers) \cite{ST1}.

\section{Holomorphic genus $2$ fibrations}

We are interested in the properties of certain specific factorizations in
$\Map_2$.

\begin{definition} Let $W_0=(T)^2$,
$W_1=(\zeta_1\cdot\zeta_2\cdot\zeta_3\cdot\zeta_4\cdot\zeta_5)^6$, and
$W_2=\sigma\cdot(\zeta_3\cdot \zeta_4\cdot \zeta_5\cdot\zeta_2\cdot
\zeta_3\cdot\zeta_4\cdot\zeta_1\cdot\zeta_2\cdot\zeta_3)^2\cdot (T)$,
where $T=\zeta_1\cdot\zeta_2\cdot\zeta_3\cdot\zeta_4\cdot\zeta_5\cdot
\zeta_5\cdot \zeta_4\cdot \zeta_3\cdot\zeta_2\cdot\zeta_1$.
\end{definition}

In this definition the notation $(\cdots)^n$ means that the sequence of
Dehn twists is repeated $n$ times. It is fairly easy to check that $W_0$,
$W_1$ and $W_2$ are all factorizations of the identity element of $\Map_2$
as a product of 20, 30, and 29 positive Dehn twists respectively (for
$W_0$ and $W_1$ this follows immediately from the presentation of $\Map_2$;
see below for $W_2$). 

\begin{figure}[t]
\begin{center}
\setlength{\unitlength}{1cm}
\begin{tabular}{ccc}
\begin{picture}(3.5,2)(0,-1)
\put(0,-1){\line(1,0){3.5}}
\put(0,1){\line(1,0){3.5}}
\put(0,-1){\line(0,1){2}}
\put(3.5,-1){\line(0,1){2}}
\multiput(0.25,-0.62)(0,0.25){6}{\line(1,0){3}}
\put(2,-0.9){\line(0,1){1.8}}
\put(2.5,-0.9){\line(0,1){1.8}}
\end{picture}
&
\begin{picture}(3.6,2)(0,-1)
\put(0,-1){\line(1,0){3.6}}
\put(1,1){\line(1,0){2.6}}
\put(0,-1){\line(1,2){1}}
\put(3.6,-1){\line(0,1){2}}
\put(1,0.5){\line(5,-2){2.4}}
\put(1,0.3){\line(4,-1){2.4}}
\put(1,-0.05){\line(1,0){2.4}}
\put(1,-0.53){\line(5,2){2.4}}
\put(1,-0.27){\line(4,1){2.4}}
\put(1,-0.8){\line(3,2){2.4}}
\end{picture}
&
\begin{picture}(3.6,2)(0,-1)
\put(0,-1){\line(1,0){3.6}}
\put(1,1){\line(1,0){2.6}}
\put(0,-1){\line(1,2){1}}
\put(3.6,-1){\line(0,1){2}}
\put(1,0.46){\line(5,-2){2.4}}
\put(1,0.433){\line(3,-1){2.4}}
\put(1,0.4){\line(4,-1){2.4}}
\put(1,-0.46){\line(5,2){2.4}}
\put(1,-0.433){\line(3,1){2.4}}
\put(1,-0.4){\line(4,1){2.4}}
\put(3,-0.9){\line(0,1){1.8}}
\put(1.4,-0.3){\circle*{0.06}}
\put(1.4,0.3){\circle*{0.06}}
\end{picture}
\\
$D_0\subset \F_0$
&
$D_1\subset \F_1$
&
$D_2\subset \F_1$
\end{tabular}
\caption{}
\end{center}
\end{figure}

\begin{lemma}
The factorization $W_0$ describes the genus 2 Lefschetz fibration $f_0$ on
the rational surface obtained as a double covering of 
$\F_0=\CP^1\times\CP^1$ branched along a smooth algebraic curve $B_0$ of
bidegree $(6,2)$. The factorization $W_1$ corresponds to the genus 2
Lefschetz fibration $f_1$ on the blown-up K3 surface obtained as a
double covering of $\F_1=\CP^2\#\overline{\CP} {}^2$ branched along a smooth
algebraic curve $B_1$ in the linear system $|6L|$, where $L$ is a line
in $\CP^2$ avoiding the blown-up point.
\end{lemma}

\proof $B_0$ can be degenerated into a singular 
curve $D_0$ consisting of 6 sections and 2 fibers intersecting in 12
nodes (see Figure 2). We can recover $B_0$ from $D_0$ by first smoothing
the intersections of the first section with the two fibers, giving us a
component of bidegree $(1,2)$, and then smoothing the remaining 10 nodes,
each of which produces two vertical tangencies. The braid factorization
corresponding to $B_0$ can therefore be expressed as
$((x_1)^2 \cdot (x_2)_{x_1}^2\cdot (x_3)_{x_2x_1}^2\cdot
(x_4)_{x_3x_2x_1}^2\cdot (x_5)_{x_4x_3x_2x_1}^2)^2$, or equivalently
after suitable Hurwitz moves, $(x_1\cdot x_2\cdot x_3\cdot x_4\cdot
x_5\cdot x_5\cdot x_4\cdot x_3\cdot x_2\cdot x_1)^2$. Lifting this
braid factorization to the mapping class group, we obtain $W_0$ as claimed.
Alternately, it is easy to check that the braid factorization for a smooth
curve of bidegree $(6,1)$ is $x_1\cdot x_2\cdot x_3\cdot x_4\cdot
x_5\cdot x_5\cdot x_4\cdot x_3\cdot x_2\cdot x_1$;
we can then conclude by observing that $B_0$ is 
the fiber sum of two such curves.

In the case of the curve $B_1$, by definition the braid monodromy is
exactly that of a smooth plane curve of degree $6$ as defined by Moishezon
in \cite{Moi}; it can be computed e.g.\ by degenerating $B_1$ to a union
of 6 lines in generic position ($D_1$ in Figure 2), and is known to be given
by the factorization $(x_1\cdot x_2\cdot x_3\cdot x_4\cdot x_5)^6$. Lifting
to $\Map_2$, we obtain that the monodromy factorization
for the corresponding double branched covering is exactly $W_1$.
\endproof

\begin{lemma} Let $\tau\in\Map_2$ be a Dehn twist, and let $F$ be a 
factorization of a central element of $\Map_2$. If $F\sim \tau\cdot
F'$ for some $F'$, then the factorization $(F)_\tau$ 
obtained from $F$ by simultaneous conjugation of all factors by $\tau$
is Hurwitz equivalent to $F$.
\end{lemma}

\proof We have: 
$(F)_\tau \sim \tau\cdot (F')_\tau\sim F'\cdot \tau\sim (\tau)_{F'}
\cdot F'=\tau\cdot F'\sim F$. The first and last steps follow from the
assumption; the second step corresponds to moving $\tau$ to the right
across all the factors of $(F')_\tau$, while in the third step all the 
factors of $F'$ are moved to the right across $\tau$. Also observe that
$(\tau)_{F'}=\tau$ because the product of all factors in $F'$ commutes
with $\tau$. \endproof

\begin{lemma} The factorizations $T$, $W_0$ and $W_1$ are {\em fully
invariant}, i.e.\ for any element $\gamma\in\Map_2$ we have 
$(T)_\gamma\sim T$, $(W_0)_\gamma\sim W_0$, and
$(W_1)_\gamma\sim W_1$.
\end{lemma}

\proof It is obviously sufficient to prove that $(T)_{\zeta_i}\sim T$
and $(W_1)_{\zeta_i}\sim W_1$ for all $1\le i\le 5$. By moving
the first $\zeta_i$ factor in $T$ or $W_1$ to the left, we obtain a
Hurwitz equivalent factorization of the form $\zeta_i\cdot \ldots$;
therefore the result follows immediately from Lemma 6.
\endproof

A direct consequence of Lemma 7 is that all fiber sums of the holomorphic
fibrations $f_0$ and $f_1$ (with monodromies $W_0$ and $W_1$) are 
{\it untwisted}. More precisely, when two Lefschetz fibrations with
monodromy factorizations $F$ and $F'$ are glued to each other along a fiber,
the resulting fibration normally depends on the isotopy class $\phi$ of a
diffeomorphism between the two fibers to be identified, and its monodromy
is given by a factorization of the form $(F)\cdot (F')_\phi$. However,
when the building blocks are made of copies of $f_0$ and $f_1$, Lemma 7 
implies that the result of the fiber sum operation is independent of the 
chosen identification diffeomorphisms; e.g., we can always take $\phi=1$.

\begin{lemma}
$(W_1)^2\sim (W_0)^3$.
\end{lemma}

\proof Let $\rho=\zeta_1\zeta_2\zeta_3\zeta_4\zeta_5
\zeta_1\zeta_2\zeta_3\zeta_4\zeta_1\zeta_2\zeta_3\zeta_1\zeta_2\zeta_1$
be the reflection of the genus 2 surface $\Sigma$ about its central axis.
It follows from Lemma 7 that 
$(W_1)^2\sim W_1\cdot (W_1)_{\rho}=(\zeta_1\cdot\zeta_2\cdot
\zeta_3\cdot\zeta_4\cdot \zeta_5)^6\cdot (\zeta_5\cdot\zeta_4\cdot\zeta_3
\cdot \zeta_2\cdot \zeta_1)^6$. The central part of this factorization is
exactly $T$; after moving it to the right, we obtain the new identity
$(W_1)^2\sim (\zeta_1\cdot\zeta_2\cdot
\zeta_3\cdot\zeta_4\cdot \zeta_5)^5\cdot (\zeta_5\cdot\zeta_4\cdot\zeta_3
\cdot \zeta_2\cdot \zeta_1)^5\cdot 
(T)_{(\zeta_5\zeta_4\zeta_3\zeta_2\zeta_1)^5}$. Repeating the same
operation four more times, we get
$(W_1)^2\sim \prod_{j=0}^5 (T)_{(\zeta_5\zeta_4\zeta_3\zeta_2\zeta_1)^j}$.
Using the full invariance property of $T$ (Lemma 7), it follows that
$(W_1)^2\sim (T)^6=(W_0)^3$.

A more geometric argument is as follows: $(W_1)^2$ is the monodromy
factorization of the fiber sum $f_1\# f_1$ of two copies of $f_1$, 
which is a double covering of the fiber sum of two copies of $(\F_1,B_1)$.
Therefore $f_1\# f_1$ is a double covering of $(\F_2,B')$ where
$\F_2=\mathbb{P}(O\oplus O(2))$ is the second Hirzebruch surface and $B'$ is
a smooth algebraic curve in the linear system $|6S|$, where $S$ is a section
of $\F_2$ ($S\cdot S=2$). On the other hand $(W_0)^3$ is the monodromy
factorization of the fiber sum $f_0\# f_0\# f_0$, which is a double covering
of the fiber sum of three copies of $(\F_0,B_0)$, i.e.\ a double covering of
$(\F_0,B'')$ where $B''$ is a smooth algebraic curve of bidegree $(6,6)$.
The conclusion follows from the fact that $(\F_2,B')$ and $(\F_0,B'')$ are
deformation equivalent.
\endproof

We can now reformulate the holomorphicity result obtained by Siebert and 
Tian~\cite{ST2} in terms of mapping class group factorizations. Say that a
factorization is {\it transitive} if the images of the factors under the
morphism $\Map_2\to S_6$ mapping $\zeta_i$ to the transposition $(i,i+1)$ 
generate the entire symmetric group $S_6$.

\begin{theorem}[Siebert-Tian \cite{ST2}]
Any transitive factorization of the identity element as a product of
positive Dehn twists along non-separating curves in $\Map_2$ is Hurwitz
equivalent to a factorization of the form $(W_0)^k\cdot (W_1)^\epsilon$
for some integers $k\ge 0$ and $\epsilon\in \{0,1\}$.
\end{theorem}

What Siebert and Tian have shown is in fact that any such factorization
is the monodromy of a holomorphic Lefschetz fibration, which can be realized
as a double covering of a ruled surface branched along a smooth connected
holomorphic curve intersecting the generic fiber in $6$ points. However,
we can always assume that the ruled surface is either $\F_0$ or $\F_1$
(either by the topological classification of ruled surfaces or using
Lemma 8). In the first case, the branch curve has bidegree $(6,2k)$ for 
some integer $k$, and the corresponding monodromy is $(W_0)^k$, while 
in the second case the branch curve realizes the homology class
$6[L]+2k[F]$ for some integer $k$ (here $F$ is a fiber of $\F_1$), and
the corresponding monodromy is $(W_0)^k\cdot W_1$.

We now look at examples of genus 2 fibrations with reducible singular
fibers.

\begin{definition}
Let $B_2\subset \F_1$ be an algebraic curve in the linear system
$|6L+F|$, presenting two triple points in the same fiber $F_0$.
Let $P_2$ be the surface obtained by blowing up $\F_1$ at the two triple 
points of $B_2$, and denote by $\hat{B}_2$ and $\hat{F}_0$ the proper 
transforms of $B_2$ and $F_0$ in $P_2$. Consider the double covering 
$\pi:\hat{X}_2\to P_2$ branched along $\hat{B}_2\cup \hat{F}_0$, and let 
$X_2$ be the surface obtained by blowing down the $-1$-curve
$\pi^{-1}(\hat{F}_0)$ in $\hat{X}_2$.
\end{definition}

Let us check that this construction is well-defined. The easiest
way to construct the curve $B_2$ is to start with a curve $C$ of degree 7 in
$\CP^2$ with two triple points $p_1$ and $p_2$. If we choose $C$
generically, we can assume that the three branches of $C$ through $p_i$
intersect each other transversely and are transverse to the line $L_0$
through $p_1$ and $p_2$. Therefore the line $L_0$ intersects $C$ 
transversely in another point $p$, and by blowing up $\CP^2$ at $p$ we obtain 
the desired curve $B_2$ (see also below). Next, we blow up the two triple
points $p_1$ and $p_2$, which turns $B_2$ into a smooth curve $\hat{B}_2$,
disjoint from $\hat{F}_0$. Denoting by $E_1$ and $E_2$ the exceptional
divisors of the two blow-ups, we have $[\hat{B}_2]=6[L]+[F]-3[E_1]-3[E_2]$ 
and $[\hat{F}_0]=[F]-[E_1]-[E_2]$, so that $[\hat{B}_2]+[\hat{F}_0]=
6[L]+2[F]-4[E_1]-4[E_2]$ is divisible by $2$; therefore the double covering
$\pi:\hat{X}_2\to P_2$ is well-defined.

The complex surface $\hat{X}_2$ is equipped with a natural holomorphic genus 2
fibration $\hat{f}_2$, obtained by composing $\pi:\hat{X}_2\to P_2$ with 
the natural projections to $\F_1$ and then to $S^2$. The fiber of
$\hat{f}_2$ corresponding to $F_0\subset\F_1$ consists of three components:
two elliptic curves of square $-2$ obtained as double coverings of the
exceptional curves $E_1$ and $E_2$ in $P_2$, with 4 branch points in each
case (three on $\hat{B}_2$ and one on $\hat{F}_0$), and 
a rational curve of square $-1$, the preimage of $\hat{F}_0$.
After blowing down the rational component, we obtain on $X_2$ a holomorphic
genus 2 fibration $f_2$, with one reducible fiber consisting of two elliptic
components. It is easy to check that near this singular point $f_2$ presents
the local model expected of a Lefschetz fibration, and that the vanishing
cycle for this fiber is the loop obtained by lifting any simple closed loop
that separates the two triple points of $B_2$ inside the fiber $F_0$ of
$\F_1$.

\begin{lemma}
The complex surface $X_2$ carries a natural holomorphic genus 2 Lefschetz
fibration, with monodromy factorization $W_2$.
\end{lemma}

\proof
We need to calculate the braid monodromy factorization associated to the
curve $B_2\subset \F_1$. For this purpose, observe that $B_2$ can be
degenerated to a union $D_2$ of 6 lines in two groups of three, $L_1,L_2,L_3$
and $L_4,L_5,L_6$, and a fiber $F$, with two triple points and 15 nodes 
(cf.\ Figure 2). The monodromy around the fiber containing the two
triple points is given by the braid $\delta=(x_1x_2)^3(x_4x_5)^3$. The
9 nodes corresponding to the mutual intersections of the two groups of three
lines give rise to 18 vertical tangencies in $B_2$, and the corresponding
factorization is $$x_3^2\cdot (x_4)_{x_3}^2\cdot (x_5)_{x_4x_3}^2\cdot
[x_2^2\cdot (x_3)_{x_2}^2\cdot (x_4)_{x_3x_2}^2]_{(x_5x_4x_3)}\cdot
[x_1^2\cdot (x_2)_{x_1}^2\cdot (x_3)_{x_2x_1}^2]_{(x_4x_3x_2x_5x_4x_3)}.$$
After suitable Hurwitz moves, this expression can be rewritten as
$$x_3\cdot x_4\cdot (x_5)^2\cdot x_4\cdot x_3\cdot [x_2\cdot x_3\cdot
(x_4)^2\cdot x_3\cdot x_2]_{(x_5x_4x_3)}\cdot [x_1\cdot x_2\cdot (x_3)^2
\cdot x_2\cdot x_1]_{(x_4x_3x_2x_5x_4x_3)},$$ or equivalently as
$x_3\cdot x_4\cdot x_5\cdot x_2\cdot x_3\cdot x_4\cdot x_1\cdot x_2\cdot
x_3\cdot x_3\cdot x_2\cdot x_1\cdot x_4\cdot x_3\cdot x_2 \cdot x_5\cdot x_4
\cdot x_3$, which is in turn equivalent to
$(x_3\cdot x_4\cdot x_5\cdot x_2\cdot x_3\cdot x_4\cdot x_1\cdot
x_2\cdot x_3)^2$. Finally, the six intersections of the lines
$L_1,\dots,L_6$ with the fiber $F$ give rise to 10 vertical tangencies,
for which the same argument as for Lemma 5 gives the monodromy factorization
$x_1\cdot x_2\cdot x_3\cdot x_4\cdot x_5\cdot x_5\cdot x_4\cdot x_3\cdot
x_2\cdot x_1$. We conclude by lifting the monodromy of $B_2$ to the mapping
class group, observing that the contribution $\delta$ of the fiber
containing the triple points lifts to the Dehn twist $\sigma$.
\endproof

\begin{theorem}
Fix integers $m\ge 0$, $\epsilon\in \{0,1\}$ and $k\ge \frac{3}{2}m+1$. Then the
Hirzebruch surface $\F_{m+\epsilon}=\mathbb{P}(O\oplus O(m+\epsilon))$
contains a complex curve $B_{k,\epsilon,m}$ in the linear system
$|6S+(m+2k)F|$ (where $S$ is a section of square $(m+\epsilon)$ and $F$
is a fiber), having $2m$ triple points lying in $m$ distinct fibers
of $\F_{m+\epsilon}$ as its only singularities.

Moreover, after blowing up $\F_{m+\epsilon}$ at the $2m$ triple points,
passing to a double covering, and blowing down $m$ rational $-1$-curves,
we obtain a complex surface and a holomorphic genus 2 fibration
$f_{k,\epsilon,m}:X_{k,\epsilon,m}\to S^2$ with monodromy factorization
$(W_0)^{k}\cdot (W_1)^{\epsilon}\cdot (W_2)^m$.
\end{theorem}

\proof
We first construct the curve $B_{k,\epsilon,m}$ by perturbation of a
singular configuration $D_{k,\epsilon,m}$ consisting of 6 sections of
$\F_{m+\epsilon}$ together with $m+2k$ fibers.
Since the case of smooth curves is a classical result, we can assume that
$m\ge 1$. Also observe that, since the intersection number of
$B_{k,\epsilon,m}$ with a fiber is equal to $6$, the $2m$ triple points
must come in pairs lying in the same fiber: $p_{2i-1}$,
$p_{2i}\in F_i$, $1\le i\le m$.

Let $u_0$ and $u_1$ be generic sections of the line bundle 
$O_{\CP^1}(m+\epsilon)$, without common zeroes. Define six sections $S_{\alpha,\beta}$
($\alpha\in\{0,1\}$, $\beta\in \{0,1,2\})$ of $\F_{m+\epsilon}$
as the projectivizations of the sections $(1,(-1)^\alpha u_0+c_{\alpha\beta} 
u_1)$ of $O_{\CP^1}\oplus O_{\CP^1}(m+\epsilon)$, where $c_{\alpha\beta}$ are
small generic complex numbers. The three sections $S_{0,\beta}$ intersect
each other in $(m+\epsilon)$ triple points, in the fibers
$F_1,\dots,F_{m+\epsilon}$ above the points of $\CP^1$ where $u_1$
vanishes, and similarly for the three sections $S_{1,\beta}$;
generic choices of parameters ensure that all the other intersections between
these six sections are transverse and lie in different fibers
of $\F_{m+\epsilon}$. We define $D_{k,\epsilon,m}$ to be the singular
configuration consisting of the six sections $S_{\alpha,\beta}$ together
with $m+2k$ generic fibers of $\F_{m+\epsilon}$ intersecting the sections
in six distinct points.

Let $s\in H^0(\F_{m+\epsilon},O(mF))$ be the product of the sections of
$O(F)$ defining the fibers $F_1,\dots,F_m$ containing $2m$ of the triple
points of $D_{k,\epsilon,m}$.
Let $s'$ be a generic section of the line bundle $L=O(D_{k,\epsilon,m}-4mF)=
O(6S+(2k-3m)F)$ over $\F_{m+\epsilon}$. Because $k\ge \frac{3}{2}m+1$, 
the linear system $|L|$ is base point free, and so we can assume that $s'$ 
does not vanish at any of the double or triple points of $D_{k,\epsilon,m}$.
Finally, let $s_0$ be the section defining $D_{k,\epsilon,m}$. 

We consider the section $s_\lambda=s_0+\lambda s^4 s'\in
H^0(\F_{m+\epsilon}, O(6S+(2k+m)F))$,
where $\lambda\neq 0$ is a generic small complex number.
Because the perturbation $s^4 s'$ vanishes at order $4$ at each of the
triple points $p_1,\dots,p_{2m}$ of $D_{k,\epsilon,m}$ in the
fibers $F_1,\dots,F_m$, it is easy to check that all the curves 
$D(\lambda)=s_\lambda^{-1}(0)$ present triple points at $p_1,\dots,p_{2m}$.
On the other hand, since $s^4 s'$ does not vanish at any of the other
singular points of $D_{k,\epsilon,m}$, for generic $\lambda$ the
curve $D(\lambda)$ presents no other singularities than the triple points
$p_1,\dots,p_{2m}$ (this follows e.g.\ from Bertini's theorem); this
gives us the curve $B_{k,\epsilon,m}$ with the desired properties.
Moreover, generic choices of the parameters ensure that the
vertical tangencies of $B_{k,\epsilon,m}$ all lie in distinct fibers of
$\F_{m+\epsilon}$; in that case, the double covering construction will
give rise to a Lefschetz fibration.

The braid monodromy of the curve $B_{k,\epsilon,m}$ can be computed using
the existence of a degeneration to the singular configuration
$D_{k,\epsilon,m}$ (taking $\lambda\to 0$ in the above construction); a
calculation similar to the proofs of Lemma 5 and Lemma 11 yields that it
consists of $k$ copies of the braid factorization of the curve $B_0$,
$\epsilon$ copies of the braid factorization of $B_1$, and $m$ copies of 
the braid factorization of $B_2$. Another way to see this is to observe
that the surface $\F_{m+\epsilon}$ admits a decomposition into a fiber sum
of $k$ copies of $\F_0$ and $m+\epsilon$ copies of $\F_1$, in such a way
that the singular configuration $D_{k,\epsilon,m}$ naturally decomposes
into $k$ copies of $D_0$, $\epsilon$ copies of a degenerate version of $D_1$
presenting some triple points, and $m$ copies of
$D_2$. After a suitable smoothing, we obtain that the pair
$(\F_{m+\epsilon}, B_{k,\epsilon,m})$ splits as the untwisted fiber sum
$k\,(\F_0,B_0)\# \epsilon\,(\F_1,B_1)\# m\,(\F_1,B_2)$.

By the same process as in the construction of the surface $X_2$, we can 
blow up the $2m$ triple points of $B_{k,\epsilon,m}$, take a double covering
branched along the proper transforms of $B_{k,\epsilon,m}$ and of the fibers
through the triple points, and blow down $m$ rational components, to obtain
a complex surface $X_{k,\epsilon,m}$ equipped with a holomorphic genus 2
Lefschetz fibration $f_{k,\epsilon,m}:X_{k,\epsilon,m}\to S^2$. 
Because of the structure of $B_{k,\epsilon,m}$, it is
easy to see that $f_{k,\epsilon,m}$ splits into an untwisted fiber sum of 
$k$ copies of $f_0$, $\epsilon$ copies of $f_1$, and $m$ copies of $f_2$. 
Therefore, its monodromy is described by the factorization
$(W_0)^k\cdot (W_1)^\epsilon\cdot (W_2)^m$.
\endproof

\section{Proof of the main result}

In order to prove Theorem 1, we will use the following lemma which
allows us to trade one reducible singular fiber against a collection of
irreducible singular fibers:

\begin{lemma}
$(\zeta_1\cdot\zeta_2)^3\cdot(\zeta_4\cdot\zeta_5)^3\cdot (T)\cdot (W_2)
\sim \sigma\cdot (W_0)\cdot (W_1)$.
\end{lemma}

\proof Let $\Phi=\zeta_3\cdot\zeta_4\cdot\zeta_5\cdot\zeta_2\cdot\zeta_3
\cdot\zeta_4\cdot\zeta_1\cdot\zeta_2\cdot\zeta_3$, and observe that 
$(\zeta_i)_\Phi=\zeta_{6-i}$ for $i\in\{1,2,4,5\}$. Therefore, we have
$$(\zeta_4\cdot\zeta_5)^3\cdot (\Phi)^2 \sim (\Phi)\cdot \zeta_1\cdot
\zeta_2\cdot\zeta_1\cdot (\Phi)\cdot \zeta_5\cdot\zeta_4\cdot\zeta_5\sim
(\Phi)\cdot (\zeta_1\cdot\zeta_2\cdot\zeta_3\cdot\zeta_4\cdot\zeta_5)^3.$$
Moreover, $(\zeta_1\cdot \zeta_2)^3\cdot (\Phi)\sim
\zeta_1\cdot\zeta_2\cdot\zeta_1\cdot (\Phi)\cdot \zeta_5\cdot\zeta_4\cdot
\zeta_5\sim (\zeta_1\cdot\zeta_2\cdot\zeta_3\cdot\zeta_4\cdot\zeta_5)^3.$
It follows that $(\zeta_1\cdot\zeta_2)^3\cdot (\zeta_4\cdot\zeta_5)^3
\cdot (\Phi)^2\sim W_1$. Recalling that 
$(\zeta_1\zeta_2)^3(\zeta_4\zeta_5)^3=\sigma I$ and $W_2=\sigma\cdot
(\Phi)^2\cdot (T)$, and using the invariance property of $T$ (Lemma 7), 
we have
$$(\zeta_1\cdot\zeta_2)^3\cdot(\zeta_4\cdot\zeta_5)^3\cdot (T)\cdot (W_2)
\sim \sigma\cdot (T)\cdot (\zeta_1\cdot\zeta_2)^3\cdot(\zeta_4\cdot\zeta_5)^3
\cdot (\Phi)^2 \cdot (T)\sim \sigma\cdot (T)\cdot (W_1)\cdot (T),$$
which is Hurwitz equivalent to $\sigma\cdot (T)^2\cdot (W_1)=\sigma
\cdot (W_0)\cdot (W_1)$.
\endproof

\proof[Proof of Theorem 1]
We argue by induction on the number $m$ of reducible singular fibers.
If there are no separating Dehn twists, then after summing with at least
one copy of $W_0$ to ensure transitivity, we obtain a transitive
factorization of the identity element into Dehn twists along non-separating
curves, which by Theorem 9 is of the expected form.

Assume that Theorem 1 holds for all factorizations with $m-1$ separating
Dehn twists, and consider a factorization $F$ with $m$ separating Dehn
twists. By Hurwitz moves we can bring one of the separating Dehn twists to
the right-most position in $F$ and assume that $F=(F')\cdot \tilde\sigma$,
where $\tilde\sigma$ is a Dehn twist about a loop separating two genus 1
components. Clearly, there exists an element $\phi\in\Map_2$ such that
$\tilde\sigma=(\sigma)_{\phi^{-1}}$. Using the relation
$I^2=1$, we can express each $\zeta_i^{-1}$ as a product of the generators
$\zeta_1,\dots,\zeta_5$, and therefore $\phi$ can be expressed as a positive
word involving only the generators $\zeta_1,\dots,\zeta_5$ (and not their
inverses).

Starting with the factorization $\tilde\sigma\cdot (W_0)^n$, we can 
selectively move $\sigma$ to the right across the various factors $W_0$, 
sometimes conjugating the factors of $W_0$ and sometimes conjugating
$\tilde\sigma$. If we choose the factors by which we conjugate
$\tilde\sigma$ according to the expression of $\phi$ in terms of
$\zeta_1,\dots,\zeta_5$, and if $n$ is sufficiently large, we obtain
that $\tilde\sigma\cdot (W_0)^n\sim (F'')\cdot \sigma$, for some factorization
$F''$ involving only non-separating Dehn twists. Therefore, using Lemma 8
and Lemma 13 we have
$$F\cdot (W_0)^{n+4}\sim F'\cdot \tilde\sigma\cdot (W_0)^n\cdot W_0\cdot
(W_1)^2\sim F'\cdot F''\cdot \sigma\cdot W_0\cdot (W_1)^2\sim
\tilde{F}\cdot W_2,$$ where $\tilde{F}=F'\cdot F''\cdot (\zeta_1
\cdot \zeta_2)^3\cdot (\zeta_4\cdot\zeta_5)^3\cdot T\cdot W_1$.
Next we observe that $\tilde{F}$ is a factorization of the identity element
with $m-1$ separating Dehn twists, therefore by assumption there exist
integers $\tilde{n},k,\epsilon$ such that $\tilde{F}\cdot (W_0)^{\tilde{n}}
\sim (W_0)^{k}\cdot (W_1)^\epsilon\cdot (W_2)^{m-1}$. It follows
that $F\cdot (W_0)^{n+\tilde{n}+4}\sim \tilde{F}\cdot W_2\cdot
(W_0)^{\tilde{n}}\sim \tilde{F}\cdot (W_0)^{\tilde{n}}\cdot W_2\sim
(W_0)^k \cdot (W_1)^\epsilon\cdot (W_2)^m$. This concludes the proof, since
it is clear that the splitting remains valid after adding extra copies of 
$W_0$.
\endproof

\proof[Proof of Corollary 2]
First of all, as observed at the beginning of \S 2 we can assume that $f$
is relatively minimal, i.e.\ all vanishing cycles are homotopically
non-trivial. Let $F$ be a monodromy factorization corresponding to $f$,
and observe that by Theorem 1 we have a splitting of the form
$F\cdot (W_0)^n\sim (W_0)^{n+k}\cdot (W_1)^{\epsilon}\cdot (W_2)^m$.
If $n$ is chosen large enough then $n+k\ge \frac{3}{2}m+1$, and so by
Theorem 12 the right-hand side is the monodromy of the holomorphic fibration
$f_{n+k,\epsilon,m}$, while the left-hand side corresponds to the fiber sum
$f\# n\,f_0$.
\endproof

\noindent {\bf Remark.} Pending a suitable extension of the result of
Siebert and Tian to higher genus hyperelliptic Lefschetz fibrations with
transitive monodromy and irreducible singular fibers, the techniques
described here can be generalized to higher genus hyperelliptic fibrations 
almost without modification. The main difference is the existence of different
types of reducible fibers, classified by the geni $h$ and $g-h$ of the two
components; this makes it necessary to replace $W_2$ with a larger collection
of building blocks, obtained e.g.\ from complex curves in $\F_1$ that
intersect the generic fiber in $2g+2$ points and present two multiple points
with multiplicities $2h+1$ and $2(g-h)+1$ in the same fiber.

\end{document}